\documentclass[11pt]{article}
\usepackage[paperheight=30cm,paperwidth=20cm,textwidth=15cm]{geometry}
\usepackage{amsfonts}
\usepackage{amsmath}
\usepackage{mathrsfs}
\usepackage{caption}
\usepackage{multicol}
\usepackage{blindtext}
\usepackage{multicol}
\usepackage{amssymb}
\usepackage{graphicx}
\usepackage{pdfpages}
\def\qed{\hfill $\square$ \\}

\def\Z{\mathbb Z}
\def\Aut{{\rm Aut}}

\newtheorem{theorem}{Theorem}[section]
\newtheorem{lemma}[theorem]{Lemma}%
\newtheorem{corollary}[theorem]{Corollary}%
\newtheorem{definition}[theorem]{Definition}%
\newtheorem{example}[theorem]{Example}%
\newtheorem{remark}[theorem]{Remark}%
\def\pf{\noindent{\it Proof.} }
 
\begin{document}



\title{{\bf On automorphism groups of polar codes}}


\author{Jicheng Ma \\[+3pt]
Chongqing Key Lab. of Group {\&} Graph Theories and Applications,\\ Chongqing University of Arts and Sciences,\\
Chongqing 402160, China \\[+3pt]   {\normalsize ma{\_}jicheng@hotmail.com}\\[+10pt]
Guiying Yan \\[+3pt]
Academy of Mathematics and System Science,\\ Chinese Academy of Sciences,\\  Beijing 100080, China \\[+3pt] {\normalsize yangy@amss.ac.cn}}

\date{}

\maketitle


{\bf Abstract}:
Over the past years, Polar codes have arisen as a highly effective class of linear codes, equipped with a decoding algorithm of low computational complexity. 
This family of codes share a common algebraic formalism with the well-known Reed-Muller codes, which involves monomial evaluations.  
As useful algebraic codes, more specifically known as decreasing monomial codes, a lot of decoding work has been done on Reed-Muller codes based on their rich code automorphisms. 
In 2021, a new permutation group decoder, referred to as the automorphism ensemble (AE) decoder, was introduced.  This decoder can be applied to Polar codes and has been shown to produce similar decoding effects.
However, identifying the right set of code automorphisms that enhance decoding performance for Polar codes remains a challenging task. 
This paper aims to characterize the full automorphism group of Polar codes. 
We will prove a reduction theorem that effectively reduces the problem of determining the full automorphism group of arbitrary random Polar codes to that of a specified class of Polar codes.
Besides, we give exact classification of the full automorphism groups of families of Polar codes that are constructed using the Reed-Muller codes.


 \section{Introduction}
 
 Polar codes, first introduced by Ar{\i}kan \cite{Arikan}, are used as a very powerful family of codes in error correction and source coding. It is a class of linear block codes relying on the phenomenon of channel polarization. 
 In particular,  they are shown to be capacity-achieving on binary memoryless symmetric channels under successive cancellation (SC) decoding for infinite block length. However, in the finite-length performance, SC decoding is not performed as expected.  To address this issue, the successive cancellation list (SCL) decoding was proposed in \cite{TalVardy},  which enhances the performance of SC decoding under finite-length conditions.
Furthermore, cyclic redundancy check (CRC)-aided (CA-SCL) decoding scheme was proposed in \cite{NiuChen}, achieving excellent performance at block error probability for polar code decoding. 
In order to avoid the extra decoding delay due to decoding operations pausing, information exchange among parallel SC decoders in hardware implementation of CA-SCL decoders, permutation-based decoders that run independently in parallel on permuted codewords were introduced in \cite{Korada} and \cite{HKU}. This approach has been further proposed in \cite{DHMG} for SC under stage permutations on the factor graph; a similar approach was proposed in \cite{PCB} for cancellation (SCAN) and in \cite{EECB} for belief propagation (BP).

In 2021, a new and effective permutation group decoder for Reed-Muller codes \cite{Muller, Reed}, termed automorphism ensemble (AE) decoder, was introduced based on the automorphism group of the code \cite{GEECB}. In stead of only stage permutations, the AE decoder employs more SC-invariant automorphisms. 
An crucial aspect of AE decoding lies in the identification and avoidance of SC-invariant automorphisms. 
The automorphism group of Reed-Muller codes is known and isomorphic to the general affine group \cite{BDOT}, which admits a large order of automorphism group for AE decoding. 
Existing evidence suggests that the application of AE decoding to polar codes can 
give similar beneficial effects on both codes. And this for sure requires a thorough understanding of the code automorphisms of polar code. 
Until now, there has been an amount of research conducted on characterizing the code automorphisms of polar codes.  For instance, in 2016, Bardet et al. \cite{BDOT} proved that the automorphisms formed by lower-triangular affine (LTA) transformations were  subgroups of automorphism groups of polar codes. Subsequently, in 2021, Geiselhart et al. \cite{GEECB}  furthered this research by proving that these automorphisms were 
SC-invariant, indicating that SC-based AE decoding does not utilize these permutations. 
Additionally, \cite{GEECB} and \cite{LZLWTYM} respectively confirmed that the block lower-triangular affine (BLTA) group is the complete affine automorphism group of polar codes.
Furthermore, it has been verified in \cite{GEECB, PBL1, PBL2} that BLTA transformations are better performed under AE decoding. In particular, in \cite{PBL2}, Pillet et al.  classified the affine automorphism group into equivalent classes, demonstrating that two equivalent automorphisms produce identical SC decoding results.  

Previous researches have been conducted to identify SC-invariant affine automorphisms for general polar codes \cite{GEECB, PBL2}. A necessary and sufficient condition related to the block lower-triangular structure of transformation matrices was given \cite{YLZLWYM}. Additionally,  Bioglio et al. \cite{BPL} provided an analysis of the algebraic properties of the affine automorphism group of polar codes, including a novel description of its structure and classification.
 However, there are currently limited known results regarding automorphisms outside the affine automorphism group, and few methods available for selection. This is the primary motivation for this work. 
 
This paper focuses on polar code automorphisms that fall outside the affine automorphism subgroup, also known as non-affine automorphisms.  
We will prove a reduction theorem that effectively reduces the problem of determining the full automorphism group of arbitrary random polar codes to that of determining the automorphism group of a specific class of polar codes.
Additionally, we give exact classification of the full automorphism groups of families of polar codes that are constructed using the Reed-Muller codes. 
The paper is organized as follows. In Section 2, we provide a brief overview of polar codes, monomial codes, Reed-Muller codes, and code automorphism and so on. We also present several useful results that will be used in followed sections. In Section 3, we present our reduction theorem, and in Section 4 we classify the full automorphism group of a family of polar codes that is constructed using the Reed-Muller codes. 

  \section{Preliminaries}


\subsection{Monomial Codes}


Let $\{x_0, x_1, \cdots, x_{n-1}\}$ be a collection of $n$ variables of values in $\mathbb{F}_2$.  
For any binary $n$-tuple $v= (v_0, v_1, \cdots, v_{n-1})\in\mathbb{F}_2^n$, let $wt(\cdot)$ denote the (Hamming) weight. Then the expression  $$x^v= \prod_{i=0}^{n}x_{i}^{v_i}$$ denotes a monomial of degree $wt(v)$. Its evaluation vector $ev(x^v)\in \mathbb{F}_2^n$ can be obtained by evaluating $x^v$ at all points of $\mathbb{F}_2^n$. 
Monomial codes of length $N = 2^n$ are a family of linear codes that can be derived from evaluations of boolean functions in $n$ variables. In particular, we will use $\bf 0$ and $\bf 1$ to represent monomials of weight 0 and $2^n$, respectively. 
The negative boolean variable $\bar x_i$ is defined as $\bar x_i=\neg x_i=  {\bf 1}\oplus  x_i$. For the sake of simplicity throughout this paper, we will replace $\bar x_i$ with $x_i$. 
Note that for a monomial codeword $x_i$ of length $2^{n}$, the concatenation $(x_i\ |\ x_i)$ of $x_i$ with itself  is actually the monomial codeword $x_i$ of length $2^{n+1}$. In particular,  $x_{n}= ({\bf 1}\ |\ {\bf 0})$. This may give a recursive definition of monomial codes. For example, Table 1 lists all monomials of length $2^4$. 

 \begin{center}{\small 
\captionof{table}{All monomials of length $2^4$}
\begin{tabular}{ cl cl}
 \hline 
{\bf 1} & 11111111\,11111111 & $x_0$ & 10101010\,10101010 \, \\[+2pt]\hline
  $x_1$ & 11001100\,11001100 &
 $x_2$ & 11110000\,11110000 \\ \hline   $x_3$ & 11111111\,00000000 &
$x_0x_1$ & 10001000\,10001000 \\ \hline   $x_0x_2$ & 10100000\,10100000 &
$x_1x_2$ & 11000000\,11000000   \\ \hline  
$x_0x_3$ & 10101010\,00000000  &
$x_1x_3$ & 11001100\,00000000   \\ \hline   $x_2x_3$ & 11110000\,00000000 & 
$x_0x_1x_2$ & 10000000\,10000000   \\ \hline   $x_0x_1x_3$ & 10001000\,00000000 &
$x_0x_2x_3$ & 10100000\,00000000   \\ \hline   $x_1x_2x_3$ & 11000000\,00000000 &
$x_0x_1x_2x_3$ & 10000000\,00000000  \\ \hline
\end{tabular}}
 \end{center}


\medskip
A length $2^n$ monomial code, say $\mathscr M_n$,  of dimension $K$ is generated by $K$ monomials, and those $K$ chosen monomials form the {\em generating monomial set}, say $M$,  of the code, while the linear combinations of their evaluations over $\mathbb F_2$ provide the {\em codebook} of the code.  However, the family of monomial codes is too large to investigate, and thus one would include a special family of monomial codes named decreasing monomial codes as follows. 
The definition of decreasing monomial codes can be seen in lots of references, see \cite{GEECB} for example. 

\begin{definition}[Decreasing Monomial Codes]
A monomial code of length $2^n$ with generating monomial set $M$ is called {\em decreasing} if  
$\forall g\in M$, $\forall f\in {\mathscr M}_{n}$ with $f\preccurlyeq g\Rightarrow f\in M$, where the partial order $`\preccurlyeq'$ of monomials is defined as 
$x_{i_1}\cdots x_{i_s}\preccurlyeq x_{j_1}\cdots x_{j_s}\Leftrightarrow i_k\leqslant j_k$ and 
$f\preccurlyeq g\Leftrightarrow \exists g'\, |\, g$ with deg($g'$)= deg($f$) and $f\preccurlyeq g'$. 
\end{definition}

In fact, several significant monomial codes are decreasing monomial codes, and we will include two classes in the following context. 

\subsection{Reed-Muller codes}

The Reed-Muller (RM) codes are one of the oldest families of algebraic codes named after D. E. Muller \cite{Muller} who introduced the codes in 1954, and I. S. Reed \cite{Reed} who proposed the first efficient decoding algorithm. The RM codes can be described in several different but equivalent ways, see \cite{MacWilliamsSloane} for example, and here, we refer the monomial codes way. 

The RM codes  of length $2^n$, denoted by  ${\rm RM}(r, n)$,  are monomial codes generated by all monomials up to degree $r$. Clearly, the RM codes are decreasing monomial codes. 
For instance, the ${\rm RM}(1, n)$ are monomial codes with generating monomial set $M= \{ {\bf 1}, x_0, \cdots, x_{n-1}\}$.  
%

\subsection{Polar Codes}
A polar code of length $N= 2^n$ and dimension $K$ is a binary linear block code defined by a transformation matrix $T_N= T_2^{\otimes n}$, where the kernel matrix $T_2\triangleq \tiny\begin{pmatrix} 1 &0 \\1 &1\end{pmatrix}$ and ${\otimes n}$ denotes $n$ times Kronecker product of $T_2$ with itself,  
an information set $\mathcal I\in \Z_{N}$ and a frozen set $\mathcal F= \Z_{N}\setminus\mathcal I$ where $\Z_{N}= \{0, 1, \cdots, N-1\}$, that is an additive group of integers modulo $N$. The encoding, an input vector ${\bf u}= (u_0, u_1, \cdots, u_{N-1})$ is defined by assigning $u_i= 0$ for $i\in \mathcal F$, which will be called frozen bits, and memorizing information in the other non-frozen bits from $\mathcal I$. Then the code is given by ${\mathscr C}= \{ x= u\cdot T_N\, |\, u\in {\mathbb F}_2^{N}\}$. The information set $\mathcal I$ and frozen set $\mathcal F$ are usually selected according to the reliabilities of the virtual bit-channels resulting from the polarization. 

It can be seen that each row of $T_N$ is in fact the evaluation vector of monomial $x^v$. Hence, polar codes are actually monomial codes. 
More specifically, 
in \cite{BDOT}, Bardet et al. showed that Polar codes are decreasing monomial codes. 
\begin{theorem}\cite[Theorem 1]{BDOT} Polar codes are decreasing monomial codes. 
\end{theorem}


\subsection{Code automorphisms}
A {\em permutation} $\pi$ on finite set $\Omega= \{1, 2, \cdots, N\}$ is a bijection from $\Omega$ to itself, and the {\em symmetric group}, say ${\rm S}_{N}$,  is the set of all permutations. 
For a linear block code $\mathscr{C}$ of length $N= 2^n$, an {\em automorphism} $\tau$ of $\mathscr{C}$ is a permutation on the $2^n$ codeword bits that maps any codeword $x\in \mathscr{C}$ into another codeword $x'=\tau(x)\in \mathscr{C}$. The full {\em automorphism group} $\Aut{\mathscr{C}}$ of $\mathscr{C}$ is the set containing all automorphisms of code $\mathscr{C}$ using composition as the multiplication. 
Hence, for a monomial code $\mathscr{M}$ and each automorphism $\tau\in \Aut\mathscr{M}$, one can see that 

$$\tau (x_i+x_j)= \tau (x_i)+\tau (x_j)$$
and
$$\tau (x_ix_j)= \tau (x_i) \tau (x_j).$$


In the following, we will consider several known code automorphisms of some classical decreasing monomial codes. First of all, it is known that the RM codes are affine-invariant in the sense that the automorphism group contains a subgroup isomorphic to the affine linear group ${\rm AGL}(n, 2)$, more specifically, we have 

\begin{theorem}\cite[Ch. 13. Theorem 24]{MacWilliamsSloane} 

(a) $\Aut({\rm RM}(r, n))= {\rm AGL}(n, 2)$ for $1\leqslant r\leqslant n-2$, otherwise 

(b) $\Aut({\rm RM}(r, n))= {\rm S}_{2^{n}}$, the symmetric group of degree $2^n$. 
\end{theorem}


In fact, for each affine automorphism $\tau\in {\rm AGL}(n, 2)$, it can be defined as the  linear transformation of $n$ variables by 

$$\tau:\ \  x\rightarrow x'= Ax+b$$ where $x, x'\in \mathbb{F}_2^n$, and $A$ is an $n\times n$ binary invertible matrix and $b$ is a binary column vector of length $n$. Equivalently, we have 

$$\tau:\ \ \begin{pmatrix} {\bf 1}\\ x \end{pmatrix}\rightarrow \begin{pmatrix} {\bf 1}\\x'\end{pmatrix}=  \begin{pmatrix} 1&0\\b&A\end{pmatrix}\begin{pmatrix} {\bf 1}\\ x\end{pmatrix}.$$ Hence, this gives a $(n+1)$ by $(n+1)$ matrix presentation of $\tau$.

 \medskip
In \cite{BDOT}, Bardet et al. mentioned the following remark, which states a theoretical motivation of this work.  

\begin{remark} Although the automorphism group of Reed-Muller codes is well-known, the question remains open for decreasing monomial codes.
\end{remark}

 \medskip
 The following example shows an easy attempt of checking and calculating code automorphisms. 
  \begin{example}\label{exam} 
Let $\mathscr M$ be a monomial code of length $8$ with generating monomial set $M= \{{\bf 1}, x_0, \cdots, x_{m}\}$ where $m\le 2$. 
Let $\tau= (1, 2, 3, \cdots, 8)\in {\rm S}_{8}$ be a permutation on the eight codeword bits, and $\tau^2= (1, 3, 5, 7)(2, 4, 6, 8)$. Then, one can see that the actions of $\tau$ and $\tau^2$ on monomials are as follows. 
  
$${\tau}(x_0)= {\tau}(1010 1010)= (0101 0101)= {\bf 1}+ x_0, $$
$${\tau}(x_1)= {\tau}(1100 1100)= (1001 1001)= {\bf 1}+x_0+x_1,$$
and thus 
${\tau^2}(x_1)=  {\bf 1}+x_1$.  
Hence, when $m= 1$, the permutation $\tau$ is an automorphism of $\mathscr{M}$. 
Similarly, by calculation we have  
\begin{align*}
{\tau}(x_2)&= \ 
{\bf 1}+x_0+x_1+ x_0x_1+x_2\\
{\tau^2}(x_2)&= \  
 {\bf 1}+x_1+x_2. \end{align*}
Hence, when  $m= 2$, the permutation $\tau$ is a code automorphism but not affine  of $\mathscr{M}$. In addition, $\tau^2$ is an affine automorphism. Furthermore, we have 

$$\tau^2:\ \ \begin{pmatrix} {\bf 1}\\x_0\\ x_1\\ x_2\end{pmatrix}\rightarrow \begin{pmatrix} {\bf 1}\\x_0\\ {\bf 1}+x_1\\ {\bf 1}+x_1+x_2\end{pmatrix}=   \begin{pmatrix}
1&0&0&0\\0&1 &0 & 0 \\ 
1&0& 1 &0  \\  
1&0 & 1 & 1    
\end{pmatrix}\begin{pmatrix}
{\bf 1}\\
x_0  \\ 
x_1 \\  
x_2    
\end{pmatrix}.$$


Furthermore, when $m= 1$. 
Let $\eta\in  \Aut\mathscr M$ be an automorphism which satisfies that 
\begin{align*}
{\eta}(x_0)&= \ 
{\bf 1}+x_0+x_1\\
{\eta}(x_1)&= \  
 {\bf 1}+x_1. \end{align*}
Additionally, suppose that $\eta(x_2)= {\bf 1}+x_0+x_1+ x_0x_1+x_2$.  The method of determining permutation presentation can be described as follows. 
%
%
In $\eta(x_0)$, one can see that the first bit is of value 1. Thus 
$\eta$ maps the 1st, 3rd, 5th  or 7th bit, which are of value 1 in $x_0$, to the 1st bit. 
In the meantime, in $\eta(x_1)$, one can see that the first bit is of value 0. Thus $\eta$ maps the 3rd or 7th bit, which are of value 0 in $x_1$, to the 1st bit. 
\begin{center}{\small 
\begin{tabular}{ c | c  c c c}
  & images under the action of $\eta$    \\ \hline
$x_0= (1010 1010)$ & $(1001 1001)= {\bf 1}+x_0+x_1$    \\  [+5pt]
$x_1= (11001100)$ &$(0011 0011)= {\bf 1}+x_1$    \\    [+5pt]
$ x_2= (11110000)$ &$  (11100001)= {\bf 1}+x_0+x_1+ x_0x_1+x_2$     \\  \hline
\end{tabular}}
\end{center} 
In addition, in $\eta(x_2)$, one can see that the first bit is of value 1. Hence $\eta$ maps the 3rd bit, which is of value 1 in $x_{2}$, to the first bit.  
Similarly, one can work out the exact permutation presentation of $\eta$ which is $\eta= (1, 8, 6, 7, 5, 4, 2, 3)$. 
\end{example}

    
\section{A reduction theorem}
 
Let ${\mathcal P_{m} }$ be any polar code of length $2^m$ with 
  $$P_m=\{ {\bf 1}, x_0, x_1, \cdots, x_{m-1}\}\cup S(x_0, x_1, \cdots, x_{m-1})  $$ 
being the generating monomial set  where $S(x_0, x_1, \cdots, x_{m-1})$ contains certain monomials of degree greater than 1.  
 For $0\leqslant i\leqslant m-1$, by the concatenation definition of monomial codes, $(x_i\,|\, x_i)$ gives a monomial codeword  of length $2^{m+1}$. 
Now, for $n>m$,  let ${\mathcal P_n }$ be a corresponding polar code of ${\cal P}_m$ of length $2^n$ with generating monomial set 
$$P_n=\{ {\bf 1}, x_0', x_1', \cdots, x_{m-1}'\}\cup S(x_0', x_1', \cdots, x_{m-1}') $$ where $x_i'= (x_i\,|\,x_i \, \cdots\, x_i\, |\, x_i)$, a concatenation of $2^{n-m}$ number of $x_i\in P_m$. 
By taking each monomial codeword $x_i'$ of $\mathcal P_n$ as a $2^{n-m}\times 2^{m}$ array as follows 
 
$$x_i' =\begin{pmatrix}
x_i\\
x_i\\
\cdots \\
x_i\\
x_i  
\end{pmatrix}\ {\rm where}\  0\leqslant i\leqslant m-1. $$ 
Then one can verify the following theorem. 
\begin{theorem}\label{r} The full automorphism group $\Aut\mathcal P_n$ of $\mathcal P_n$ is isomorphic to $({\rm S}_{2^{n-m}})^{2^m}: G$ where $G\cong \Aut\mathcal P_{m}$, the full automorphism group of $\mathcal P_m$. 
\end{theorem} 
\pf  The $2^{n-m}\times 2^{m}$ array presentation of each codeword of $\mathcal P_n$ reveals that all rows are identical. Thus, each row shift of a column results in a code automorphism, which generates a subgroup isomorphic to the direct product of $2^m$ number of symmetric groups of degree $2^{n-m}$, say  $({\rm S}_{2^{n-m}})^{2^m}$. In particular, each element in $({\rm S}_{2^{n-m}})^{2^m}$ fixes the monomials $x_i'$ where $0\leqslant i\leqslant m-1$. 

Besides, we claim that all code automorphisms which acting as  shifts of columns generate a subgroup isomorphic to $\Aut\mathcal P_m$. 
Let $\rho_{i}$ be a mapping which maps each coordinate $j$ to $j+(i-1)\cdot2^{m}$ with $1\leqslant j\leqslant 2^m$ and $1\leqslant i\leqslant n-m$.
For each permutation $\pi\in {\rm S}_{2^m}$, let $\tilde{\pi}\in {\rm S}_{2^n}$ be a lifting permutation of $\pi$ such that $\tilde{\pi}= \rho_1(\pi)\rho_2(\pi)\cdots\rho_{n-m}(\pi)$. 
Then it is easy to see  that $\tilde{\pi}\in \Aut\mathcal P_n$ if and only if $\pi\in \Aut\mathcal P_m$. 
Let $G$ be the subgroup generated by all such $\tilde{\pi}\in \Aut\mathcal P_n$, and thus $G\cong \Aut\mathcal P_m$. Moreover, one can see that the automorphism subgroup $({\rm S}_{2^{n-m}})^{2^m}$ of row shifts can be normalized by $G$. Thus $\mathcal P_n$ contains a subgroup $({\rm S}_{2^{n-m}})^{2^m}: G$ of automorphisms. 
In the meantime, suppose $\rho_i(\pi)\in  \Aut\mathcal P_n$. Then $\rho_i(\pi)$ has to fix each codeword  since $\rho_i(\pi)$ only permutes the columns of a fixing row $i$. Hence, $\pi$ has to fix each codeword of $\mathcal P_m$, which suggests that  both $\pi$ and $\rho_i(\pi)$ are the trivial permutation. The claim holds. 

In the following, we will show that there exists no other code automorphisms. On the contrary, consider $\varphi\in {\rm S}_{2^{n}}\setminus (({\rm S}_{2^{n-m}})^{2^m}:G)$, and suppose $\varphi\in \Aut\mathcal P_n$. Then $\varphi$ has to permutes columns of non-identical rows. Again, as all of the $2^{n-m}$ rows of each codeword of $\mathcal P_n$ are identical, thus for each $x'\in \mathcal P_n$ such that $${x'}=\begin{pmatrix}
x\\
x\\
\cdots \\
x\\
x  
\end{pmatrix},$$  one has 

$${x'}^{\varphi} =\begin{pmatrix}
x\\
x\\
\cdots \\
x\\
x  
\end{pmatrix}^{\varphi}= \begin{pmatrix}
y\\
y\\
\cdots \\
y\\
y  
\end{pmatrix}\in \mathcal P_n.$$ 
Note that, by the assumption $x, y\in{\mathcal P}_m$, there must exist an automorphism $g\in \Aut{\mathcal P}_m$ such that $x^g=y$.
Hence, similar to the above argument, there exists a lifting automorphism $\tilde{g}\in G$ of $g$ such that ${x'}^{\varphi}= {x'}^{\tilde{g}}$. 
Moreover, there must exist an automorphism $\tau\in ({\rm S}_{2^{n-m}})^{2^m}$ such that ${\varphi}= \tilde{g}\tau\in ({\rm S}_{2^{n-m}})^{2^m}: G$, a contradiction. 
Therefore, one has $\Aut\mathcal P_n\cong ({\rm S}_{2^{n-m}})^{2^m}: G$.  
 \qed

Especially, we have the followed remark.

 \begin{remark} Let $\mathcal P_m$ be the Reed-Muller code ${\rm RM}(r, m)$, by Theorem~3.1, the full automorphism group of ${\cal P}_n$ is one of the following:
 
 (a) $\Aut({\cal P}_n)\cong ({\rm S}_{2^{n-m}})^{{2^{m}}}: {\rm AGL}(m, 2)$ for $1\leqslant r\leqslant m-2$, or  

(b) $\Aut({\cal P}_n)\cong ({\rm S}_{2^{n-m}})^{{2^{m}}}:  {\rm S}_{2^{m}}$ for $r= 1$ or $ m-1$.

 \end{remark}

 \smallskip
 
Additionally, by Theorem~3.1, one can see that in order to classify the full automorphism groups of the Polar codes ${\cal P}_n$,  one only need to determine the subgroup $G$. Besides, for the subgroup $({\rm S}_{2^{n-m}})^{2^m}\leqslant \Aut{\cal P}_n$, one can see that by the arguments of Theorem~3.1 each element  fixes all codewords of ${\cal P}_n$. However, the induced action of each automorphism in $({\rm S}_{2^{n-m}})^{2^m}$ might be not affine. In the following, we aim to characterize the non-affine subgroup of automorphisms in $\Aut\mathcal P_n\cong ({\rm S}_{2^{n-m}})^{2^m}: G$. 

Now assuming that $n= m+1$, then by Theorem~\ref{r}, the full automorphism group of ${\cal P}_n$ is isomorphic to $(\Z_2)^{2^{n-1}}: G$. Additionally, 
we have the following results:

\begin{lemma}\label{r-lem} When $n= m+1$, for each automorphism $\tau\in  \Aut\mathcal P_n$, the action of $\tau$ on $x_{n-1}$ gives that $\tau(x_{n-1})= x_{n-1}+c$ where $c\in \mathcal P_n$. 
\end{lemma} 
 \pf First of all, for each automorphism $\tau\in \Aut{\cal P}_n$ of ${\cal P}_n$, the image of $\tau(x_{n-1})$ has to contain a term which is divisible by $x_{n-1}$, otherwise, $x_{n-1}\ne  \tau^{-1}(\tau(x_{n-1}))$, a contradiction, as $\tau^{-1}\in \Aut{\cal P}_n$ as well.  
Suppose  that $\tau(x_{n-1})= \alpha x_{n-1}+\beta$ where $\alpha, \beta\in {\cal P}_n$. Then one can verify that $\alpha= \bf 1$. On the contrary, suppose that $\alpha$ is of degree greater than 0.

$$x_{n-1}= \tau ^{-1}(\alpha x_{n-1})+ \tau ^{-1}(\beta) = \tau ^{-1}(\alpha)\tau ^{-1}( x_{n-1})+ \tau ^{-1}(\beta).$$ 
Note that neither $\tau ^{-1}(\alpha)$ nor $\tau^{-1}(\beta)$ contain a term of $x_{n-1}$. Thus, $x_{n-1}+ \tau ^{-1}(\beta)\ne   \tau ^{-1}(\alpha)$. 
 In the following, we will try to find a contradiction. We claim that, there exists a codeword bit of $x_{n-1}+ \tau ^{-1}(\beta)$ of value equals to 1 but at the same bit in $\tau ^{-1}(\alpha)$ is of value 0.
This is true, as the codeword $x_{n-1}+ \tau ^{-1}(\beta)$ is of weight exactly $2^{n-1}$ number of 1s, and $\tau ^{-1}(\alpha)$ contains at most $2^{n-1}$ number of 1s, moreover, in $x_{n-1}+ \tau ^{-1}(\beta)$, the values of each pair of bits 1 and $1+2^{n-1}$ are different, however, the values are the same in $\tau ^{-1}(\alpha)$. Thus the claim is true, and hence $\alpha=\bf 1$. This finishes the lemma. \qed

\begin{theorem}\label{main-1} When $n= m+1$, for each $\tau\in \Aut\mathcal{P}_n$, one has $\tau$ fixes $x_{n-1}$ or   \[\tau (x_{n-1})= x_{n-1} +  \sum_{(l_0, l_1, \cdots, l_{n-2})}^{} x_0^{l_0}x_1^{l_1}\cdots x_{n-2}^{l_{n-2}}\] where $(l_0, l_1, \cdots, l_{n-2})\in {\mathbb F}_2^{\,n}$. 
\end{theorem}
\pf  By Lemma~\ref{r-lem}, we know that for each $\tau\in \Aut\mathcal{P}_n$, $\tau (x_{n-1})= x_{n-1}+ c$ where $c\in {\cal P}_n$. 
Now, suppose that  \[ c=  \sum_{(l_0, l_1, \cdots, l_{n-2})}^{} x_0^{l_0}x_1^{l_1}\cdots x_{n-2}^{l_{n-2}} \] with $l_i\in {\mathbb F}_2$.  
In order to show that such permutation action of $\tau$ on $x_{n-1}$ still 
gives an automorphism of $\mathcal{P}_n$. Equivalently, by the definition of code automorphisms one only needs to show that the action $\tau$ on $x_{n-1}$ is a permutation. In fact, this is true for each sequence $(l_0, l_1, \cdots, l_{n-2})$ as the weight of $x_{n-1}+c$ still equals to $2^{n-1}$. 
\qed

It is worth noting that similar to the calculation of Example 2.5 and together with Theorem~\ref{main-1}, one can work out the exact permutation action of each automorphism $\tau\in \Aut(\mathcal{P}_n)$. 
With regard to the affine properties of each automorphisms, one can classify them from the perspective of theoretical derivations as follows.
 
Let  $A_f({\cal P})\leqslant  \Aut{\mathcal P}$ be the complete affine subgroup of automorphisms of Polar codes $\cal P$. Then  we have 

\begin{corollary}\label{order} The complete affine automorphism subgroup $A_f({\cal P}_n)$ is isomorphic to ${\mathbb{Z}_2}^{n}: H$ where $H\cong A_f({\cal P}_m)$. 
\end{corollary}
\pf By Theorem~\ref{main-1}, for each automorphism $\tau\in (\Z_2)^{2^{n-1}}$, $\tau$ fixes monomials $x_0, x_1, \cdots$, $x_{n-2}$. Suppose that $\tau$ is an affine automorphism, then the action of $\tau$ on $x_{n-1}$ must be affine. Hence, by Theorem~\ref{main-1}, $\tau$ generates a subgroup isomorphic to the elementary abelian 2-group ${\Z_2}^{n}$. Thus the corollary holds.  
\qed

 In the following, we re-define the polar codes ${{\cal P}_{m}}$ by letting 
the generating monomial set $P_m$ to be $$P_m=\{ {\bf 1}, x_0, x_1, \cdots, x_{m-1}\}\cup S(x_0, x_1, \cdots, x_{k-1})  $$ where $k<m$. 
Moreover, we would also assume that $x_{0}x_{k-1}\in S(x_0, x_1, \cdots, x_{k-1})$, and thus $x_0x_{i}\in S(x_0, x_1, \cdots, x_{k-1})$ where $1\leqslant i\leqslant k-2$.  
 Let ${\cal P}_k$ be a sub-code of ${\cal P}_{m}$ with generating monomial set 
 $P_k= P_m\setminus\{x_{k}, \cdots, x_{m-1}\}$ and is of dimension $K$.  
  Additionally, let $\bar{{\cal P}}_k$ be the corresponding Polar codes of length $2^{k}$ of ${\cal P}_k$ and with the same generating monomial sets.  
Then the following theorem holds.

 \begin{theorem} For $k\geqslant 2$, the full automorphism group $\Aut\mathcal P_{m}$ of $\mathcal P_{m}$ is a subgroup of $\Aut\mathcal P_k$, the full automorphism group of $\mathcal P_k$, and  $\Aut{\cal P}_{m}\cong (\Z_2)^{{(m-k)K}}: ({\rm GL}(m-k, 2)\times T)$ where $T\cong \Aut{\bar {\cal P}_{k}}$, the full automorphism group of ${\bar {\cal P}_{k}}$. 
\end{theorem} 
 \pf   
Let $\alpha\in \Aut{\cal P}_{m}$ be an automorphism of ${\cal P}_{m}$, we claim that for any $x_i\in {\cal P}_{m}$ with $0\leqslant i\leqslant k-1$, $\alpha(x_i)$ contains no terms of $x_s$ with $k\leqslant s\leqslant m-1$. 
First of all, we will show that $\alpha(x_i)$ contains no terms of $x_{m-1}$. On the contrary, suppose there exists an $x_i$ such that 
$\alpha(x_i) = x_{m-1}+\beta_1$ where $\beta_1$ contains no terms of $x_{m-1}$,  and suppose that $\alpha(x_0)= \beta_2$. Note that, clearly $\beta_2$ cannot equal to $\bf 1$, $\bf 0$ or $x_{m-1}+\beta_1$. 
Then $\alpha(x_0x_i)= \beta_2x_{m-1}+\beta_1\beta_2$. As $\alpha(x_0x_i)\in {\cal P}_{m}$, thus $\beta_2= x_{m-1}$ and it follows that $\beta_1= \bf 0$, a contradiction. 
 Similar to the above argument, together with the fact that ${\cal P}_m$ contains no monomials of degree greater than 1 which can be divisible by $x_{k}, x_{k+1}, \cdots, x_{m-1}$, one can verify the claim. 
 
 By Theorem~\ref{r}, the automorphism group of $\Aut\mathcal P_{k}$ is isomorphic to $({\rm S}_{2^{m-k}})^{2^{k}}: T$ where $T$ is isomorphic to the full automorphism group of the corresponding Polar codes $\bar {\cal P}_{k}$ of length $2^{k}$. 
 In the following, we are going to show which of the automorphisms of $\mathcal P_{k}$ are automorphisms of $\mathcal P_{m}$.

First of all, by the arguments of Theorem~\ref{r}, one can see that the subgroup $({\rm S}_{2^{m-k}})^{2^{k}}\leqslant \Aut\mathcal P_{k}$ fixes each monomial $x_i$ with $0\leqslant i\leqslant k-1$. 
 Hence, we only need to distinguish all those permutations that map $x_k, x_{k+1}, \cdots, x_{m-1}$ to an element inside ${\cal P}_m$.  However, this could be quite complicated. 
 For example, the permutation $\tau= (2^{m-2}+1,2^{m-1}+1 )\in ({\rm S}_{2^{m-k}})^{2^{k}}$ satisfies that $$({x_{m-1}})^{\tau}=\begin{pmatrix}
1\, 1\, \cdots 1\, 1\\
\ \ \cdots \\
1\, 1\, \cdots 1\, 1\\
0\, 0\, \cdots 0\, 0\\
\ \ \cdots\\
0\, 0\, \cdots 0\, 0
\end{pmatrix}^{\tau} = \begin{pmatrix}
1\, 1\, \cdots 1\, 1\\
\ \ \cdots \\
0\, 1\, \cdots 1\, 1\\
1\, 0\, \cdots 0\, 0\\
\ \ \cdots\\
0\, 0\, \cdots 0\, 0
\end{pmatrix} \notin\Aut{\cal P}_{m-1}. $$
Instead, now, for each automorphism $\tau\in \Aut{\cal P}_m\cap ({\rm S}_{2^{m-k}})^{2^{k}}$, suppose that 
$$(x_{j})^{\tau}= a_{j,k}x_{k}+a_{j,k+1}x_{k+1}+\cdots+a_{j,m-1}x_{m-1}+ \beta_j$$ with $k\leqslant j\leqslant m-1$, $\beta_j\in {\cal P}_k$ and $a_{j, s}\in \mathbb{F}_2$ where $k\leqslant s\leqslant m-1$.
For those $m-k$ number of vectors
\begin{center}
 $(a_{k,k}, a_{k,k+1},\cdots, a_{k, m-1})$
 
$(a_{k+1,k}, a_{k+1,k+1},\cdots,  a_{k+1, m-1})$\\
$\cdots$\\

$(a_{m-1,k}, a_{m-1,k+1},\cdots, a_{m-1, m-1})$
\end{center}
 we claim that they are linearly independent. Otherwise 
 $$(x_k+x_{k+1}+\cdots+x_{m-1})^{\tau}= \beta_k+\beta_{k+1}+\cdots+\beta_{m-1},$$ and thus 
 $$ (\beta_k+\beta_{k+1}+\cdots+\beta_{m-1})^{{\tau^{-1}}}= x_k+x_{k+1}+\cdots+x_{m-1}, $$ a contradiction. 
 Hence, those $m-k$ number of vectors generate a linear group isomorphic to ${\rm GL}(m-k, 2)$. 
 And for $k\leqslant j\leqslant m-1$, by letting $\beta_j=\bf 0$, one can see that $T\times {\rm GL}(m-k, 2)$ is a subgroup of $\Aut{\cal P}_{m}$.  
 
Furthermore, for any $\beta_j\in {\cal P}_k$ and $a_{j,k}x_{k}+a_{j,k+1}x_{k+1}+\cdots+a_{j,m-1}x_{m-1}$ as above, 
one can see that the weight of $a_{j,k}x_{k}+a_{j,k+1}x_{k+1}+\cdots+a_{j,m-1}x_{m-1}+ \beta_j$ equals to the weight of $x_j$. 
Hence, with the actions of $\tau$ on $x_0, x_1, \cdots, x_{m-1}$, the exact permutation presentation of $\tau$ can be workout. Clearly, different assignment of $\beta_j$ will result in different $\tau$.  
Now let $N\leqslant \Aut{\mathcal P}_m$ be a subgroup consists of each automorphism $\tau\in\Aut{\mathcal P}_m$ satisfies that ${x_i}^\tau= x_i$ when $0\leqslant i\leqslant k-1$ and ${x_i}^\tau= x_i+\beta_i$ with $\beta_i\in{\mathcal P}_k$ when $k\leqslant i\leqslant m-1$. Then for each $0\leqslant j\leqslant m-1$, $(x_j)^{\tau^2}= x_j$. And for $\tau_1, \tau_2\in N$, one can see that $\tau_1\tau_2= \tau_2\tau_1$, which suggests that $N$ is an elementary abelian 2-group of order $2^{(m-k)\cdot K}$. 
Let $g_1g_2\in G\times {\rm GL}(m-k, 2)$ that for $0\leqslant i\leqslant k-1$, ${x_i}^{g_1}\in {\cal P}_k$ and ${x_i}^{g_2}= x_i$, and for $k\leqslant i\leqslant m-1$, ${x_i}^{g_1}= x_i$ and ${x_i}^{g_2}=  \sum a_{j,t}x_t$, where $k\leqslant t\leqslant m-1$ and $a_{j,t}\in {\mathbb F}_2$. 
Then one has $g_1g_2$ normalizes $N$, and thus  $\Aut({\cal P}_{m})= (\Z_2)^{{(m-k)K}}: ({\rm GL}(m-k, 2)\times \Aut({\bar {\cal P}_{k}}))$. This finishes the proof.  
\qed


    \begin{corollary} With the assumption of Theorem~3.6, the Polar codes $\mathcal{P}_m$ admits non-affine automorphisms. 
 \end{corollary}
 

By the assumption of Theorems~3.1 and~3.6, one can see that in order to classify the full automorphism groups of the Polar codes ${\cal P}_n$ or ${\cal P}_m$,  one only need to determine the subgroup $T$ which is isomorphic to the automorphism group of the corresponding Polar codes $\bar{\cal P}_k$ of ${\cal P}_k$. 
 Hence, from this point onward,  we will assume ${\cal P}_r$ to be a Polar code of length $2^r$ with generating monomial set 
 $$P_r= \{ {\bf 1}, x_0, x_1, \cdots, x_{r-1}\}\cup S(x_0, x_1, \cdots, x_r)$$
where $S(x_0, x_1, \cdots, x_r)$ contains monomial $x_0x_{r-1}$ and other monomials of degree greater than 1.

  \section{Polar codes constructed from Reed-Muller codes}
 
  Before starting, 
let  ${\cal M}_r$ be the Reed-Muller codes ${\rm RM}(r, n)$ with $M_r$ as the generating monomial set.  
Let ${\cal P}_r$ be Polar codes constructed from the ${\rm RM}(r, n)$, and let $P_r$ be the generating monomial set of ${\cal P}_r$. 

Now, for $1\leqslant r<n-1$, let ${\cal P}_r$ be Polar codes with generating monomial set $P_r= {M}_r\cup \{x_0x_1\cdots x_{r}\}$. 
Note that, when $r= n-1$, ${\cal P}_{n-1}$ is actually the RM codes ${\rm RM}(r, n)$. Moreover, for $r= 1$, it can be seen that by Theorem~3.6, 
the automorphism group $\Aut{\cal P}_r$ is  isomorphic to ${\Z_2}^{4(n-2)}: ({\rm S}_4\times {\rm GL}(n-2, 2))$. For the other values of $r$, we have the following results. 
 

 \begin{theorem}{\label 1} 
For $2\leqslant r\leqslant n-3$, the full automorphism group $\Aut\mathcal{P}_r$ is a subgroup of $ \Aut\mathcal{M}_r$, and thus affine. In particular, $\Aut\mathcal{P}_r$ is isomorphic to $(\Z_2)^{(n-m)(m+1)}: ({\rm AGL}(m, 2)\times GL(n-m, 2))$ where $m= r+1$.
\end{theorem}
\begin{pf}
We claim that for each $g\in \Aut\mathcal{P}_r$, and for $0\leqslant i\leqslant n-1$, $g(x_i)$ is of degree 1, and for $0\leqslant j\leqslant r$, $g(x_j)$ contains no terms of  $x_k$ where $r+1\leqslant k\leqslant n-1$. 
As $x_0x_1\cdots x_r\in \mathcal{P}_r$ which is the only monomial of degree $r+1$, thus $g(x_0x_1\cdots x_r)\in \mathcal{P}_r$ and is of degree $r+1$. Hence, the claim holds, otherwise contradicts to the above facts.  

Hence, one can see that each automorphism $g$ must be an affine automorphism, and a $n+1$ by $n+1$ matrix presentation of $g$ can be given as follow. 

$$T_g=   \begin{pmatrix}
1&0&0\\ 
b_1 &T_1& 0\\ 
b_2 &B&T_2\\ 
\end{pmatrix}. $$
where $T_1$ and $T_2$ are respectively $(r+1)$ by $(r+1)$  and $(n-r-1)$ by $(n-r-1)$ invertible matrices over $\mathbb{F}_2$, and $b_1, b_2, B$ are respectively $(r+1)$ by 1, $(n-r-1)$  by 1, $(n-r-1)$ by $(r+1)$ matrices over $\mathbb{F}_2$. 

Now, by letting $T_2$ be the identity matrix, and $b_2, B$ are  zero matrices. Then all $T_g$ generate a subgroup isomorphic to the affine linear group ${\rm AGL}(r+1, 2)$. Similarly, by letting $T_1$ be the identity matrix, and $b_1, b_2, B$ are zero matrices. Then all $T_g$ generate a subgroup isomorphic to the linear group ${\rm GL}(n-r-1, 2)$. Furthermore, by letting $T_1$ and $T_2$ be the identity matrices, and $b_1, b_2$ be the zero matrix. Then all $T_g$ generate a subgroup isomorphic to $(\Z_2)^{(n-r-1)(r+1)}$. By letting $T_1$ and $T_2$ be the identity matrices, and $b_1, B$ be the zero matrix. Then all $T_g$ generate a subgroup isomorphic to $(\Z_2)^{n-r-1}$. 
Note that, each of those five matrices are relatively independent. 
Hence, let $m= r+1$, the full automorphism group $\Aut\mathcal{P}_r$ is of order $2^{(n-m)(m+1)}\cdot |{\rm AGL}(m, 2)|\cdot |{\rm GL}(n-m, 2)|$. Moreover, one can see that $(\Z_2)^{(n-m)(m+1)}$ can be normalized by ${\rm AGL}(m, 2))\times {\rm GL}(n-m, 2)$, and thus $\Aut\mathcal{P}_r\cong (\Z_2)^{(n-m)(m+1)}:({\rm AGL}(m, 2)\times {\rm GL}(n-m, 2))$. 
\qed
 \end{pf} 

   \begin{remark}Compare to Theorem~3.6, one can see that, the automorphism groups in Theorem~3.6 and Theorem~4.2 are isomorphic. In fact, the two Polar codes are complements.
\end{remark}

Especially, when $r= n-2$, we have 
 
\begin{remark} Let $\tau\in {\rm S}_{2^{n}}$ be a permutation such that $\tau(x_{n-1})= x_{n-1}+ x_{n-3}x_{n-2}$ and $\tau(x_i)= x_i$ for $0\leqslant i\leqslant n-2$. 
Then one can easily verify that $\tau$ is an automorphism of ${\cal P}_{n-2}$. 
Therefore, the full automorphism group $\Aut{\cal P}_{n-2}$ admits non-affine automorphisms. 

  \end{remark}

 \bigskip\bigskip
\noindent {\Large\bf Acknowledgments}

\smallskip
\noindent
The authors acknowledge the use of the {\sc Magma} computational package \cite{Magma}, which helped show the way to many of the results given in this paper.

\end{document}